\documentclass[12pt]{amsart}

\usepackage{euscript,amssymb,amsmath}
\usepackage{graphics}
\usepackage{amscd}

\advance\hoffset by -.5truein

\def\wt{\operatorname{wt}}
\def\pre{\operatorname{pre}}
\def\Hom{\operatorname{Hom}}
\def\Brace{\operatorname{Brace}}
\def\Ker{\operatorname{Ker}}
\def\gr{\operatorname{gr}}

\def\As{\mathrm{As}}
\def\Dend{\mathrm{Dend}}
\def\Vect{\mathrm{Vec}}
\def\Lie{\mathrm{Lie}}
\def\Perm{\mathrm{Perm}}
\def\Com{\mathrm{Com}}

\newtheorem{thm}{Theorem}
\newtheorem{lem}{Lemma}
\newtheorem{rem}{Remark}

\newtheorem{cor}{Corollary}

\theoremstyle{definition}
\newtheorem{defn}{Definition}
\newtheorem{exmp}{Example}

\title{Gr\"obner---Shirshov bases for pre-associative algebras}
\author{Pavel Kolesnikov}
\address{Sobolev Institute of Mathematics, Novosibirsk, Russia}

\begin{document}

\maketitle

\section{Introduction}

The method of Gr\"obner bases developed by B. Buchberger in \cite{Buch70} is known to be one of 
the most effective techniques for 
solving various problems in different areas of mathematics wherever 
one needs to decide whether a given polynomial belongs to an ideal 
generated by a given family of polynomials (i.e., to solve the word problem). Here by ``polynomial'' we mean 
an element of the free associative and commutative algebra. 
Independently, A. I. Shirshov in \cite{Shir62} proposed similar technique 
for Lie algebras which may be easily reduced to 
the case of associative non-commutative algebras \cite{Bok76}. 
In contrast to the commutative Gr\"obner bases technique, the general one (Gr\"obner---Shirshov bases technique) does not provide an algorithm 
(in the classical ``finite'' sense) for solving the word problem, however, 
it is a powerful theoretic tool for studying a wide class of algebraic systems, 
see the reviews \cite{BC2014, BC2015} and references therein. 

One of the most important applications of the Gr\"obner---Shirshov bases technique 
is the proof of PBW-like  (Poincar\'e---Birkhoff--Witt) theorems for left adjoined 
functors to multiplication changing functors between varieties of linear algebras. 
Namely, suppose $\mathcal V$ and $\mathcal W$  are two (linear) operads governing 
the varieties of linear algebras called $\mathcal V$-algebras and $\mathcal W$-algebras, respectively. 
A morphism of operads $\omega : \mathcal W\to \mathcal V$ induces a functor from the variety of 
$\mathcal V$-algebras to the variety of $\mathcal W$-algebras. 
Such functors are called multiplication changing ones. 
Every multiplication changing functor has left adjoined functor which 
sends an arbitrary $\mathcal W$-algebra $A$ to its universal enveloping $\mathcal V$-algebra $U_\omega (A)$.
There exists canonical linear map $\iota : A\to U_\omega (A)$ which may not be injective, in general, such that $U_\omega (A)$
is generated by $\iota(A)$. Moreover, $U_\omega (A)$ carries a natural ascending filtration and its associated graded algebra 
$\gr U_\omega (A)$ is also a $\mathcal V$-algebra. As it was proposed in \cite{MikhShest}, let us say 
the triple $(\mathcal V, \mathcal W, \omega )$ to have {\em PBW-property} if 
\[
 \gr U_{\omega }(A) \simeq U_\omega (A^{(0)})
\]
as $\mathcal V$-algebras, where $A^{(0)}$ stands for the $\mathcal W$-algebra on same underlying 
space as $A$ with trivial (zero) operations.

It was shown in \cite{MikhShest} that many combinatorial properties of free $\mathcal V$- and $\mathcal W$-algebras 
are closely related provided that $(\mathcal V, \mathcal W, \omega )$ has the PBW-property. 
Gr\"obner---Shirshov bases technique is the traditional tool to decide whether $(\mathcal V, \mathcal W, \omega )$ has the PBW-property. 
To apply this technique, one has to describe the free $\mathcal V$-algebra, prove a version of Composition-Diamond Lemma for $\mathcal V$-algebras, 
and find Gr\"obner---Shirshov basis (GSB) of $U_\omega (A)$ for a generic $\mathcal W$-algebra $A$.

Note that $U_\omega (A)$ has a standard presentation by generators and defining relations. 
Suppose the operad $\mathcal W$ is generated by operations $\Omega $, and $\nu (f)$ stands for the arity of an operation $f\in \Omega $.
Let $X$ be a linear basis of $A$, then 
$U_\omega (A)$ may be presented as $\mathcal V\langle X\rangle /I( \omega(f)(x_1,\dots, x_{\nu(f)}) - f(x_1,\dots , x_{\nu(f)}) , x_i\in X, f\in \Omega  )$, 
where $\mathcal V\langle X\rangle $ is the free $\mathcal V$-algebra generated by $X$ and $I(S)$ denotes the ideal of $\mathcal V\langle X\rangle $
generated by $S$.
The set $S$ of defining relations stated here may not be a GSB, but if we can embed it into a GSB $\bar S$ generating the same ideal $I(S)=I(\bar S)$
then it remains to consider principal parts of the GSB $\bar S$ obtained.
If neither of these principal parts belong to $X$ and the list of all these parts does not depend on the particular multiplication table of $A$ 
then $(\mathcal V,\mathcal W, \omega )$ has the PBW-property \cite{MikhShest}.

Let us mention, for example, the paper \cite{BCY2010_LSA} 
where the Gr\"obner---Shirshov bases method was developed for right-symmetric (pre-Lie) algebras known also as pre-Lie algebras. 
As an application, it was actually shown that $(\pre\Lie, \Lie, ({-}) )$ has the PBW-property.
Here $\pre\Lie$ is the operad governing the variety of algebras with one binary multiplication such that 
\[
 (xy)z - x(yz) - (xz)y + x(zy) = 0,
\]
$\Lie$ is the operad of Lie algebras, and $\omega =({-}) $ turns $xy$ into $xy-yx$.
It turns out that, in these settings, the initial set of defining relations of the 
universal enveloping $\pre\Lie$-algebra of an arbitrary Lie algebra is a GSB.

In this paper, we consider similar problem in the associative settings and, as an application, prove that 
$(\Dend, \As, ({*}) )$ has the PBW-property, where 
$\Dend $ is the operad governing the variety of dendriform algebras \cite{Loday01}, 
i.e., algebras with two binary operations $\prec $ and $\succ$ satisfying the identities
\[
\begin{gathered}
 x\succ (y\succ z) = (x\succ y)\succ z + (x\prec y)\succ z, \\
 x\succ (y\prec z) = (x\succ y)\prec z , \\
 x\prec (y\succ z)+x\prec (y\prec z) = (x\prec y)\prec z,
\end{gathered}
\]
$\As $ is the operad of associative algebras, and $({*})$ turns $xy$ into $x\succ y + x\prec y$.

We develop Gr\"obner---Shirshov bases method for dendriform algebras and, as an application, compute the GSB of the 
universal enveloping dendriform algebra of an arbitrary associative algebra. In contrast 
to the pre-Lie case \cite{BCY2010_LSA}, multiplication table of $A$ itself is not a GSB, but we 
still present the complete system of relations.

In fact, as we explain in Section \ref{sec:PreAlgebras}, the relation between $\Lie $ and $\pre\Lie $
is very similar to the relation between $\As$ and $\Dend $, as well as between the operads governing 
Poisson and pre-Poisson algebras \cite{Aguiar00}, Jordan and pre-Jordan algebras \cite{Hou2013}. 
This is why we prefer the term ``pre-associative algebra'' instead of ``dendriform algebra'' for $\Dend$-algebras. 
\section{Gr\"obner---Shirshov bases in free algebras}\label{sec:Section1}

In this section, we recall the notion of a Gr\"obner---Shirshov basis in 
the free non-associative algebra $M\langle X\rangle $ and in the free
associative algebra $\As\langle U \rangle $. In this section, we will
follow \cite{Bok76} and \cite{BokChenHuang13} with minor changes.
We will also introduce some notations to be used later.

Let $X$ be a nonempty set equipped with a well order $\le $, and let $X^{**}$
be the set of all (nonempty) non-associative words in the alphabet $X$. 
For $u\in X^{**}$, denote by $|u|$ the length of $u$.
Define the {\em weight} $\wt(u)$ of $u\in X^{**}$ as follows: for $u=x\in X$ put $\wt(u)=(1,x)$, 
for $u=u_1u_2$ put $\wt(u) = (|u|, u_2, u_1)$. 
 Extend the initial order on $X$ to the order 
$\le $ on $X^{**}$ by induction on the length:
\[
 u\le v \iff \wt(u)\le \wt(v)\ \text{lexicographically}.
\]
Obviously, this is a {\em monomial} order, i.e., 
\[
 u\le v \Rightarrow wu\le wv,\ uw\le vw
\]
for all $u,v,w\in X^{**}$.

\begin{rem}
This definition of order on $X^{**}$ is slightly different from 
much more general one in \cite{BokChenHuang13}, where the weight of $u=u_1u_2$
is defined as $(|u|,u_1,u_2)$, but all statements remain valid.
\end{rem}

The set $X^{**} $ is a linear basis of the free non-associative algebra $M\langle X\rangle $ generated 
by $X$. Suppose $0\ne f,g\in M\langle X\rangle $, are non-associative polynomials, 
and let $\bar f$, $\bar g$ be principal parts of $f$, $g$. Assume $f$ and $g$ are monic, i.e., the coefficients 
at their principal parts are equal to~$1$.
If $\bar f = (u_1\dots u_k \bar g v_1 \dots v_l)$ with some bracketing, $u_i,v_j\in X^{**}$, then 
\[
 h = f - (u_1\dots u_k g v_1 \dots v_l) 
\]
is called the {\em composition of inclusion} relative to $w=\bar f$.

A set $S$ of monic non-associative polynomials is called a {\em Gr\"obner---Shirshov basis} in $M\langle X\rangle $
if for all $f,g\in S$ and for every their composition of inclusion $h$ relative to a word $w$ we have
\[
 h = \sum\limits_i  \alpha _i (u_{i1}\dots u_{ik_i} s_i v_{i1}\dots v_{il_i})_i,  
\]
where $\alpha_i\in \Bbbk$, $u_{ij},v_{ij}\in X^{**}$ ($k_i,l_i\ge 0$), $s_i\in S$, 
$(\dots )_i$ denote some bracketings, 
and 
$(u_{i1}\dots u_{ik_i} \bar s_i v_{i1}\dots v_{il_i})_i<w$.
This property of $h$ is denoted $h\equiv 0\pmod {S,w}$.

If a non-associative word $u\in X^{**}$ has no composition of inclusion with any of $g\in S$ then $u$ is said to be an $S$-{\em reduced} word. 

\begin{thm}[\cite{BokChenHuang13}]\label{thm:CD_nonass}
Let $S$ be a set of monic non-associative polynomials in $M\langle X\rangle $, and let $I_\cdot (S)$ be 
the ideal of $M\langle X\rangle $ generated by~$S$. 
Then the following statements are equivalent:
\begin{itemize}
 \item $S$ is a Gr\"obner---Shirshov basis in $M\langle X\rangle $;
 \item $0\ne f \in I_\cdot (S)$ implies $\bar f = (u_1\dots u_k s v_1 \dots v_l)$ for an appropriate $s\in S$;
 \item The set of all $S$-reduced words forms a linear basis of $M\langle X\mid S\rangle = M\langle X\rangle /I_\cdot (S)$.
\end{itemize}
\end{thm}

Introduce the following relation on $M\langle X\rangle $: 
\[
 f\to_S g \ \iff\ g = f-\alpha (u_1\dots u_k s v_1\dots v_l),
\]
where $\bar f = (u_1\dots u_k \bar s v_1\dots v_l)$, $s\in S$, $u_i,v_i\in X^{**}$, 
and $\alpha \in \Bbbk $ is the coefficient at $\bar f$ in $f$.

If there exists a chain $f\to_S g_1 \to_S \dots \to_S 0$ 
then we say that $f$ may be reduced to zero modulo $S$.
Denote it also by $f\to_S 0$.

\begin{cor}[Elimination of leading words]\label{cor:ELW_nonass}
Suppose $S$ is a GSB in $M\langle X\rangle $, $f\in M\langle X\rangle $.
Then $f\in I_\cdot (S)$ if and only if $f$ may be reduced to zero 
modulo $S$.
\end{cor}

Now recall the associative case \cite{Bok76}.
Given a nonempty set $U$, denote by $U^*$ the set of all associative words 
in $U$ relative to associative binary operation denoted by~$*$,
$U^*$ is a linear basis of the free associative algebra $\As\langle U \rangle $ generated by~$U$.
Denote by $|w|_*$ the length of a word $w\in U^*$, and let $U^\#$ be the set $U^*$ together with empty word $\epsilon$; assume $|\epsilon |_*=0$.
Suppose $U^*$ is equipped with a well monomial order relative to the operation~$*$ (e.g., the deg-lex order based on a well order on the 
alphabet $U$).

Let $f,g\in \As\langle U\rangle  $ be monic polynomials.
Define the following two types of compositions.
\begin{itemize}
\item[(AC1)] 
If $\bar f= w_1* \bar g * w_2$, $w_i\in U^\#$, then 
\[
h = f - w_1*g*w_2
\]
is called a {\em composition of inclusion} relative to $w=\bar f$.

\item[(AC2)]
If $\bar f =u *w_1$, $\bar g =w_2*u$, 
 $u,w_1,w_2\in U^*$,
then 
\[
 w_2*f - g*w_1 \in \As\langle U\rangle
\]
is called a {\em composition of intersection} relative to $w=w_1*\bar g = \bar f*w_2$.
\end{itemize}

A set $\Sigma \subset \As\langle U \rangle $ of monic polynomials is a Gr\"obner---Shirshov basis 
in $\As\langle U\rangle $ if for all $f,g\in \Sigma $ 
and for every their composition of inclusion $h$ relative to a word $w$ we have
\[
 h = \sum\limits_i \alpha_i w_i*s_i*w'_i, 
\]
where $\alpha_i\in \Bbbk $, $w_i,w'_i\in U^\#$, $s_i\in \Sigma $, and $w_i*\bar s_i*\bar w'_i<w$.
This property of $h$ is denoted $h\equiv 0\pmod {\Sigma ,*,w}$.

If a word $u\in U^*$ has no composition of inclusion with $g\in \Sigma $ then $u$ is said to be
a $\Sigma $-reduced word.

\begin{thm}[\cite{Bok76}]\label{thm:CD_assoc}
Let $\Sigma $ be a set of monic non-associative polynomials in $\As\langle U\rangle $, and let $I_* (\Sigma )$ be 
the ideal of $\As\langle U\rangle $ generated by~$\Sigma $. 
Then the following statements are equivalent:
\begin{itemize}
 \item $\Sigma $ is a Gr\"obner---Shirshov basis in $\As\langle U\rangle $;
 \item $f \in I_*(\Sigma )$ implies $\bar f = w*\bar s *w'$ for an appropriate $s\in \Sigma$, $w,w'\in U^\#$;
 \item The set all of $\Sigma $-reduced words forms a linear basis of
 $\As\langle U\mid \Sigma \rangle = \As\langle U\rangle /I_* (\Sigma )$.
\end{itemize}
\end{thm}

In an obvious way, one may define what does it mean that $f\in \As\langle U\rangle $
may be reduced to zero modulo $\Sigma $ (c.f. Corollary \ref{cor:ELW_nonass}). Denote it as $f\to_{\Sigma ,*} 0$.

\begin{cor}\label{cor:ELW_ass}
Suppose $\Sigma $ is a GSB in $\As\langle U\rangle $, $f\in \As\langle U\rangle $.
Then $f\in I_* (\Sigma )$ if and only if $f\to_{\Sigma ,*} 0 $.
\end{cor}

\section{Pre-algebras}\label{sec:PreAlgebras}
Throughout this section, by an operad $\mathcal P$ we mean a collection of linear $\Bbbk $-spaces  
$\mathcal P(n)$, $n\ge 1$, 
equipped with polylinear composition rule and linear action of the symmetric group $S_n$
satisfying the well-known axioms of associativity, identity, and equivariance (see, e.g., 
\cite{LodayVallette}). 
For example, the class $\Vect_\Bbbk $
of all linear $\Bbbk $-spaces with polylinear maps is a multi-category, but every 
single space $V\in \Vect_\Bbbk $ 
may be considered as an operad by means of $V(n)=\Hom_\Bbbk (V^{\otimes n}, V)$
with ordinary composition rule,
where symmetric group acts by permutation of arguments. 

Given an operad $\mathcal P$, a $\mathcal P$-algebra $A$ is a morphism of operads (i.e., a functor on multi-categories)
$\mathcal P\to \Vect_\Bbbk $. 
Later we will often identify a morphism $A$ and linear $\Bbbk $-space $A(\mathcal P)$ 
considered as an operad in $\Vect_\Bbbk $.

\begin{exmp}\label{exmp:Perm}
Consider the family of $\Bbbk $-spaces $\Perm(n)=\Bbbk^n$, $n\ge 1$, 
with standard bases $\big\{e^{(n)}_i \mid i=1,\dots, n \big\}$. 
Define composition rule and symmetric group actions in a natural way (\cite{Chap01}):
\[
e_i^{(n)}\big( e_{j_1}^{(m_1)}, \dots, e_{j_n}^{(m_n)} \big)
= e^{(m)}_{m_1 + \dots + m_{i-1} + j_i}, \quad
\sigma : e_i^{(n)} \mapsto e_{i\sigma}^{(n)},\ \sigma \in S_n.
\]
The operad obtained is denoted by $\Perm $. The class of $\Perm $-algebras is a variety 
of all associative algebras satisfying the following identity:
\[
 (xy)z=(yx)z.
\]
\end{exmp}

Recall the definition of what is a $\pre\mathcal P$-algebra for a given operad $\mathcal P$ \cite{GubKol_proc2013}.

Suppose $\Omega = \bigoplus\limits_{n\ge 1}\Omega(n)$ is a graded linear $\Bbbk $-space, where each 
$\Omega(n)$ is an $S_n$-module. Let $\mathcal F=\mathcal F_{\Omega }$ stand for the 
free operad generated by $\Omega $.

Denote by $\Omega^{(2)}$ the following graded linear $\Bbbk $-space:
\[
 \Omega ^{(2)} = \bigoplus\limits_{n\ge 1}\Omega(n)\otimes \Perm(n).
\]
Extend the action of $S_n$ by 
\[
 (\omega \otimes e^{(n)}_i)^\sigma = \omega ^\sigma \otimes e^{(n)}_i, \quad \omega \in \Omega (n), \ \sigma\in S_n, 
\]
and denote by $\mathcal F^{(2)}= \mathcal F_{\Omega^{(2)}}$ 
the free operad generated by $\Omega^{(2)}$.

Given an $\mathcal F^{(2)}$-algebra $A$ and a $\Perm$-algebra $P$, define a morphism 
\[
 P\boxtimes A : \mathcal F \to \Vect_\Bbbk 
\]
as follows:
\[
 (P\boxtimes A)(\mathcal F) = P\otimes A \in \Vect_\Bbbk ,
\]
\begin{equation}\label{eq:ProdExpansion}
 (P\boxtimes A)(\omega ): (p_1\otimes a_1,\dots , p_n\otimes a_n) \mapsto 
 \sum\limits_{i=1}^n P(e^{(n)}_i)(p_1,\dots, p_n)\otimes A(\omega \otimes e^{(n)}_i) (a_1,\dots, a_n), 
\end{equation}
for $\omega \in \Omega (n)$, $p_i\in P$, $a_i\in A$, $i=1,\dots, n$, $n\ge 1$.

\begin{defn}\label{defn:preP-algebra}
Assume an operad $\mathcal P$ is an image of $\mathcal F$, i.e., there exists an epimorphism 
$\tau: \mathcal F\to \mathcal P$. 

An $\mathcal F^{(2)}$-algebra $A$ is called a $\pre\mathcal P$-algebra 
if for every $\Perm $-algebra $P$ the morphism $P\boxtimes A$ is a $\mathcal P$-algebra 
(i.e., $\Ker\tau \subseteq \Ker(P\boxtimes A)$). 
\end{defn}

Obviously, it is enough to check this property for countably generated free $\Perm $-algebra $P$ only.

Let us consider one particular example in more details.

\begin{exmp}[Pre-associative algebras]
Suppose $\mathcal P=\As $ is the operad of associative algebras, here 
$\Omega = \Omega (2)$, $\dim\Omega (2)=2$, $\Omega (2)=\Bbbk \mu + \Bbbk\mu^{(12)}$, $(12)\in S_2$.
For an $\As$-algebra $A$, the image of $\mu $ is the binary operation of multiplication
$A\otimes A\to A$, $\mu^{(12)}$ is the opposite multiplication.

Then $\Omega ^{(2)}$ is a 4-dimensional space spanned by 
$\mu_i = \mu\otimes e^{(2)}_i$ and $\mu_i^{(12)}=\mu^{(12)}\otimes e^{(2)}_i$, $i=1,2$. In the usual notation, the images of 
$\mu_1$ and $\mu_2$
in an $\mathcal F^{(2)}$-algebra correspond to two binary products denoted $\succ $ and $\prec $, 
respectively. 

According to the general Definition~\ref{defn:preP-algebra}, an $\mathcal F^{(2)}$-algebra $A$ is a $\pre\As$-algebra 
if $P\boxtimes A$ is an associative algebra, i.e., 
\[
 (p_1\otimes a_1) ( (p_2\otimes a_2) (p_3\otimes a_3) )
=
 ((p_1\otimes a_1) (p_2\otimes a_2)) (p_3\otimes a_3) 
\]
for every $p_i\in P$, $a_i\in A$, where $P$ is an arbitrary $\Perm$-algebra.

The expansion \eqref{eq:ProdExpansion} for $P\boxtimes A$ in this case turns into 
\[
 (p\otimes a)(q\otimes b) = pq\otimes (a\succ b) + qp\otimes (a\prec b),
\]
so $P\boxtimes A$ is associative if and only if 
\begin{multline}\nonumber
p_1p_2p_3\otimes a_1\succ (a_2\succ a_3) + p_1p_3p_2\otimes  a_1\succ (a_2\prec a_3) \\
+ p_2p_3p_1\otimes a_1\prec (a_2\succ a_3) + p_3p_2p_1\otimes  a_1\prec (a_2\prec a_3) \\
=
p_1p_2p_3\otimes (a_1\succ a_2)\succ a_3 + p_2p_1p_3\otimes  (a_1\prec a_2)\succ a_3 \\
+ p_3p_1p_2\otimes (a_1\succ a_2)\prec a_3 + p_3p_2p_1 \otimes  (a_1\prec a_2)\prec a_3.
\end{multline}
Collecting similar terms, we obtain the following necessary and sufficient conditions for $A$ to be pre-associative:
\begin{equation}\label{eq:DendIdent}
\begin{aligned}
 & a_1\succ (a_2\succ a_3) = (a_1\succ a_2)\succ a_3 + (a_1\prec a_2)\succ a_3, \\
 & a_1\succ (a_2\prec a_3) = (a_1\succ a_2)\prec a_3 , \\
 & a_1\prec (a_2\succ a_3)+a_1\prec (a_2\prec a_3) = (a_1\prec a_2)\prec a_3,
 \end{aligned}
\end{equation}
i.e., $A$ is a dendriform algebra ($A\in \Dend $) in the sense of \cite{Loday01}.
\end{exmp}

\begin{exmp}[Pre-Lie algebras]
Let $\mathcal P=\Lie$ be the operad of Lie algebras. Then $\Omega = \Omega (2)$, $\dim\Omega(2)=1$, 
$\Omega (2)=\Bbbk \nu$, $\nu^{(12)}=-\nu $. 
The class of $\pre\Lie$-algebras determined by Definition~\ref{defn:preP-algebra} coincides 
with the variety of left-symmetric algebras (LSAs), also called pre-Lie algebras \cite{Gerst,Kozh,Vinberg}. 
\end{exmp}


\begin{rem}
 The operad $\pre\mathcal P$ governing the variety of all $\pre\mathcal P$-algebras coincides with 
 the Manin black product $\pre\Lie \bullet \mathcal P$ \cite{Vallette}.
\end{rem}

The structure of a pre-associative algebra may be equivalently described by means of two operations 
\begin{equation}\label{eqNewOperations}
x*y = x\succ y + x\prec y\quad \text{and} \quad xy = x\succ y.
\end{equation}
Let us rewrite axioms \eqref{eq:DendIdent} in terms of these operations:
\begin{equation}\label{eq:DendIdent-2}
\begin{aligned}
 & a_1* (a_2* a_3) = (a_1* a_2)* a_3 , \\
 & (a_1* a_2) a_3 = a_1 (a_2 a_3),\\
 & a_1 (a_2* a_3) = (a_1 a_2)* a_3 -(a_1,a_2,a_3) , \\
 \end{aligned}
\end{equation}
where $(x,y,z)=(xy)z-x(yz)$ is the associator.
For example, \eqref{eq:DendIdent-2} imply the following relation in $\pre\As$:
\begin{equation}\label{eq:DendIdent-21}
 a_1(a_2*a_3*a_4)
= (a_1a_2)*a_3*a_4 - (a_1,a_2,a_3)*a_4 +(a_1,a_2,a_3)a_4 -(a_1a_2)(a_3a_4).
\end{equation}

It is easy to see that the identities \eqref{eq:DendIdent-2} allow to rewrite 
every term in a pre-associative algebra generated by a set $X$ 
as a linear combination of the following monomials:
\begin{equation}\label{eq:pre-As_basis}
 u_1*u_2 *\dots * u_k,\quad k\ge 1,\ u_i\in X^{**},
\end{equation}
where $X^{**}$ denotes, as above, the set of all (nonempty) non-associative words in the alphabet $X$
relative to the second operation in \eqref{eqNewOperations}.

\begin{thm}\label{thm:PreAs_basis}
The set of all monomials \eqref{eq:pre-As_basis} is a linear basis 
of the free pre-associative algebra $\pre\As\langle X\rangle $ generated by 
$X$. Therefore, $\pre\As\langle X\rangle $ is isomorphic to $T_0(M\langle X\rangle )$
as a linear space,
where $M\langle X\rangle $ stands for the free non-associative algebra and 
$T_0$ denotes tensor algebra without identity.
\end{thm}

\begin{proof}
It is enough to consider the case $|X|=1$ since the operad $\pre\As $ is nonsymmetric.
One may easily calculate the number of different monomials \eqref{eq:pre-As_basis} of length $n$
is equal to the $n$th Catalan number $C_n$ which is known to be the dimension of 
$\Dend(n)$ \cite{Loday01}. Hence, the words \eqref{eq:pre-As_basis}
are linearly independent in $\Dend \langle X\rangle = \pre\As\langle X\rangle $.
\end{proof}

\begin{lem}\label{lem:preAs-Ideals}
Let $S$ be a subset of $\pre\As\langle X\rangle $. 
Denote by $I(S)$ the ideal of the free pre-associative algebra generated by $S$. 
Then $I(S)$ coincides with the ideal $I_*(\tilde S)$ generated in the free associative 
algebra $\As\langle U\rangle $, $U=X^{**}$, by the set $\tilde S$ of all elements 
\[
 (u_1\dots u_l s v_1\dots v_r),\quad u_i,v_j\in U,\ l,r\ge 0,\ s\in S
\]
with all possible bracketings.
\end{lem}

\begin{proof}
 Obviously, $I_*(\tilde S) \subseteq I(S)$. Conversely, 
$I_*(\tilde S)$ is an ideal in $\pre\As\langle X\rangle $ due to \eqref{eq:DendIdent-2}.
Hence, $I_*(\tilde S) \supseteq I(S)$.
\end{proof}

\section{Gr\"obner---Shirshov bases in pre-associative algebras}\label{sec:Section3}

Let $X$ be a nonempty set. Denote by $U$ the set $X^{**}$
of all non-associative words in the alphabet $X$.
By $U^*$ we denote the set of all associative words in $U$
relative to binary operation~$*$. Theorem \ref{thm:PreAs_basis} 
implies $U^*$ to be a linear basis of $\pre\As \langle X\rangle $.
Moreover, $\pre\As \langle X\rangle$ as an associative algebra 
with respect to $*$ is isomorphic to $\As\langle U\rangle $.

Assume $X$ is well-ordered, and this order is extended to $U$  
as in Section~\ref{sec:Section1}. Then extend the order to $U^*$
in a monomial way such that $|u|_*>1$ and $|v|_*=1$ imply $u>v$ for $u,v\in U^*$
(e.g., the standard deg-lex order has this property).
The latter condition guarantees that a polynomial $f\in \pre\As\langle X\rangle $
is $*$-free (i.e., belongs to $\Bbbk U$) if so is its principal part.

Consider the following types of compositions in pre-associative algebras.

\begin{enumerate}
 \item[(C1)] The composition of {\em $*$-inclusion} for $f,g\in \pre\As\langle X\rangle $ 
 is defined in the same way as the composition (AC1) in the free 
 associative algebra $\As\langle U\rangle $;
 \item[(C2)] The composition of {\em $*$-intersection} coincides with (AC2) in $\As\langle U\rangle $;
 \item[(C3)] Suppose $f,g\in \pre\As\langle X\rangle $, $|\bar g|_*=1$. Assume 
 $\bar f=v*u*v'$, $v,v'\in U^\#$, where $u \in \widehat{\{\bar g\}}$, i.e., 
 \[
  u = (w_1\dots w_k \bar g w_1' \dots w'_{k'}) \in U
 \]
 with respect to some bracketing $(\ldots )$, $w_i,w'_i\in U$, $k+k'>0$.
 Then 
 \[
  h = f - v*(w_1\dots w_k g w_1' \dots w'_{k'})*v'
 \]
 is called the composition of {\em $\succ$-inclusion} relative to $w=\bar f$;
 
 \item[(C4)] Let $f\in \pre\As\langle X\rangle $, $|\bar f|_*>1$. Then for every $v\in U$ the element 
 $  h = fv \in \pre\As\langle X\rangle $
is called the composition of {\em right multiplication} relative to $w=\overline{fv}$;

 \item[(C5)]
Let $f\in \pre\As\langle X\rangle $, $|\bar f|_*>1$. 
Then for every $u\in U^*$ the element 
$ h = u f \in \pre\As\langle X\rangle $
is called the composition of {\em left multiplication}  relative to $w=\overline{uf}$.
\end{enumerate}

Suppose $S\subseteq \pre\As\langle X\rangle $ is a set of polynomials.
Denote by $I(S)$ the ideal of $\pre\As \langle X\rangle$ generated by $S$.
Let us split $S$ as $S=S_0\dot\cup S_1$, where $S_0 = \{f \in S \mid |\bar f|_*=1\}$, 
i.e., $S_0$ consists of all $*$-free polynomials, $S_0= \Bbbk U\cap S$, $S_1=S\setminus S_0$.
Define 
\[
 \hat S_0 = \bigcup\limits_{k\ge 0} S_0^{(k)}, \quad 
 S_0^{(0)}=S_0, \ 
 S_0^{(k+1)} = \{ vf, fv \mid f\in S_0^{(k)}, v\in U \}.
\]
Denote 
\[
 \hat S = \hat S_0 \cup S_1.
\]

Note that $I_*(\hat S) \subseteq I(S)$ since the latter is equal to $I_*(\tilde S)$
by Lemma \ref{lem:preAs-Ideals}. 

A composition $h$ of type (C1)--(C3) relative to a word $w$ is said to be 
{\em trivial modulo $S$\/} if 
$h\equiv 0\pmod {\hat S,*,w}$.
Similarly, a composition $h$ of type (C4)--(C5) relative to $w$ 
is trivial modulo $S$ if 
\[
h= \sum\limits_i \alpha_i w_i*s_i*w'_i, 
\]
where $\alpha_i\in \Bbbk $, $w_i,w'_i\in U^\#$, $s_i\in \hat S$, 
and $w_i*\bar s_i*\bar w'_i\le w$.
Note that for a composition of type (C4) we should have $w_i=w'_i=\epsilon$ 
and $s_i\in \hat S_0$ since $|\bar h|_*=1$.

\begin{defn}\label{defn:GSB_preAs}
A set of monic pre-associative polynomials $S\subseteq \pre\As\langle X\rangle $ 
is a Gr\"obner---Shirshov basis (GSB) in $\pre\As\langle X\rangle $
if for all $f,g\in S$ every composition of type (C1)--(C5) is trivial modulo~$S$.
\end{defn}

Assume $S$ is a set of monic pre-associative polynomials.

\begin{lem}\label{lem:C1--C3}
The set $\hat S$ is a GSB in $\As\langle U\rangle $ if and only if 
for all $f,g\in S$ every their composition of type (C1)--(C3) is trivial modulo $S$. 
\end{lem}

\begin{proof}
If $\hat S$ is a GSB in $\As\langle U\rangle $ then compositions (C1)--(C3) of elements of $S$ 
coincide with compositions (AC1) or (AC2).

Conversely, since all compositions of type (C3) are trivial, 
$S_0$ is a GSB in $M\langle X\rangle = \Bbbk U$.

Consider $f,g\in \hat S$ as elements of $\As\langle U\rangle $, suppose they
have a composition $h$ of type (AC1) or (AC2) relative to a word $w\in U^*$.  
If $h$ is a composition of intersection (AC2) 
then both $f$, $g$ belong to $S_1$ and $h\equiv 0\mod {\hat S,*,w}$ by definition.
Assume $h$ is a composition of inclusion (AC1). If $f\in \hat S_0$, $g\in S_1$ (or converse) then 
$h$ coincides with a composition of type (C3) and thus $h\to_{\hat S} 0 $.
Consider the case when both $f$, $g$ belong to $\hat S_0$, and $h$ is a composition of type (AC1) relative to $w$. Then
\[
 \bar f = \bar g= w, \quad h = f - g\in I_\cdot (S_0), 
\]
where either $h=0$ or $\bar h<w$.
By Corollary \ref{cor:ELW_nonass}, $h\to_{S_0} 0$, i.e., $h\equiv 0\pmod {S_0,w}$. The latter trivially implies 
$h\equiv 0\pmod {\hat S,*,w}$.
\end{proof}

\begin{lem}\label{lem:C4--C5}
If for all $f\in S$ every composition of type (C4), (C5) is trivial modulo $S$ then $I_*(\hat S) = I(S)$. 
\end{lem}

\begin{proof}
Obviously, $I_*(\hat S)\subseteq I(S)$ for an arbitrary $S$. To prove the converse, it is
enough to show that the set
$I_*(S)$  is an ideal of $\pre\As\langle X\rangle $. 

It follows from \eqref{eq:DendIdent-2} that for every $s\in \hat S_0$, $w=v_1*\dots *v_k\in U^*$, $v_i\in U$,
we have $sw,ws \in I_*(\hat S_0)$, moreover, $ws\in \Bbbk \hat S_0$.  

Suppose $b$ is a generic element of  $I_*(\hat S)$, 
\[
 b = \sum\limits_i \alpha_i w_i*s_{0i}*w'_i + \sum\limits_j \beta_j v_j*s_{1j}*v_j'  ,
\]
where $w_i,w'_i,v_j,v'_j\in U^{\#}$, $s_{0i}\in \hat S_0$, $s_{1j}\in S_1$, 
$\alpha_i,\beta_j \in \Bbbk $.

It is easy to derive from \eqref{eq:DendIdent-21}
that 
$ w(w_i*s_{0i}*w_i')  \in I_*(\hat S_0)$
for all $w\in U^*$
since $(w,w_i,s_{0i}),(ww_i)(s_{0i}w_i') \in \hat S_0$.
Consider the second group of summands:
\[
 w(v_j*s_{1j}*v_j')=(w v_j)* s_{1j}* v_j'  - (w, v_j, s_{1j})* v_j'
+(w, v_j, s_{1j}) v_j' - (w v_j)(s_{1j} v_j')\in I_*(\hat S)
\]
provided that  $(w,v_j,s_{1j})\in I_*(\hat S)$. 
The latter holds since 
 $(wv_j)s_{1j}$ and $w (v_j s_{1j})=(w*v_j)s_{1j}$ are compositions of type (C5). 
 
Therefore, $wb\in I_*(\hat S)$. In a similar way, one may show $bw\in I_*(\hat S) $ for all $w\in U^*$.
This implies $I_*(\hat S)$ to be an ideal of pre-associative algebra $\pre\As\langle X\rangle $.
\end{proof}

\begin{rem}\label{rem:Ideals_to_Comp}
 If $\hat S$ is a GSB in $\As\langle U\rangle $ then the converse holds:
 $I_*(\hat S)=I(S)$ implies that all compositions of type (C4), (C5) 
 of elements of $S$ are trivial. 
\end{rem}

Indeed, if $f\in S$ then $wf, fu\in I(S)=I_*(\hat S)$ for all $w\in U^*$, $u\in U$. 
The latter means, in particular, $wf,fu\to_{\hat S,*} 0$ in $\As\langle U\rangle$, i.e., 
compositions (C4), (C5) are trivial modulo~$S$ in $\pre\As\langle X\rangle $.

\begin{thm}\label{thm:CD-preAs}
Suppose $S\subseteq  \pre\As\langle X\rangle $ is a set of monic pre-associative polynomials, and $I(S)$ is the 
ideal in $\pre\As\langle X\rangle $ generated by~$S$.
Then the following statements are equivalent:
\begin{itemize}
 \item $S$ is a GSB in $\pre\As\langle X\rangle $;
 \item $f \in I(S)$ implies $\bar f = w*\bar s *w'$ for an appropriate $s\in \hat S$, $w,w'\in U^\#$;
 \item The set all of $\hat S$-reduced words in $\As\langle U\rangle $ is a linear basis of
 $\pre\As\langle X\mid S \rangle = \pre\As\langle X\rangle /I(S)$.
\end{itemize}
\end{thm}

\begin{proof}
(i) $\Rightarrow $ (ii)\\
If $S$ is a GSB then $f\in I(S)$ implies $f\in I_*(\hat S)$ by Lemma \ref{lem:C4--C5}.
As $\hat S$ is a GSB in $\As\langle U\rangle $ by Lemma~\ref{lem:C1--C3},
$\bar f$ contains a subword $\bar s$ for an appropriate $s\in \hat S$.

(ii) $\Rightarrow $ (iii)\\
Obviously, every element of $\pre\As\langle X\mid S \rangle $ may be rewritten 
as a combination of $\hat S$-reduced words from $U^*$. If (ii) holds, these 
words are linearly independent.

(iii) $\Rightarrow $ (i)\\
By Lemma \ref{lem:preAs-Ideals}, $\pre\As\langle X\mid S\rangle $
is isomorphic to $\As\langle U\mid \tilde S\rangle $ as a linear space. 
As $\hat S\subseteq \tilde S$, there is a linear map 
$\alpha :\As\langle U\mid \hat S\rangle \to \pre\As\langle X\mid S\rangle$
such that the diagram 
\[
 \begin{CD}
\pre\As\langle X \rangle   @>\simeq >> \As\langle U\rangle \\
@VVV @VVV \\
\pre\As\langle X\mid S\rangle @<\alpha << \As\langle U\mid \hat S\rangle
 \end{CD}
\]
commutes.

In general, the set of all $\hat S$-reduced words is a complete (but not necessarily linearly independent) set in 
the space $\As\langle U\mid \hat S\rangle$. In the case when (iii) holds, the images of these words form a basis in 
$\pre\As\langle X\mid S\rangle$, so the condition (iii) of Theorem~\ref{thm:CD_assoc} holds for $\hat S$. Hence, 
$S$ satisfies the conditions of Lemma~\ref{lem:C1--C3}.
Moreover, $\alpha $ has to be a linear isomorphism, so $I_*(\hat S)=I_*(\tilde S)=I(S)$ by Lemma~\ref{lem:preAs-Ideals}.
Thus, Lemma \ref{lem:C4--C5} and Remark \ref{rem:Ideals_to_Comp} imply $S$ to be a GSB in $\pre\As\langle X\rangle$.
\end{proof}

Therefore, in order to check whether $S$ is a GSB in $\pre\As\langle X\rangle $, 
one has to compute all compositions of type (C1)--(C3) of all elements of $S$ first and show they are all trivial.
Next, apply elimination of leading monomial of $\hat S$ to all compositions of type (C4), (C5) to show they all 
reduce to zero.

\section{Applications}

Suppose $A\in \pre\As$ is defined as a quotient of $\pre\As\langle X\rangle $ relative to 
an ideal $I$. If $I=I(S)$ for a GSB  $S$ then, by abuse of terminology,
we say $S$ is a GSB of the pre-associative algebra~$A$.

In this section, we compute Gr\"obner---Shirshov bases for two series of pre-associative algebras: 
for free Zinbiel algebras and for universal pre-associative envelopes of associative algebras. 

\subsection{GSB of the free pre-commutative algebra}

Recall that Zinbiel algebra is the name proposed by J.-M. Lemaire (see \cite{Loday01}) for commutative dendriform algebras, 
i.e., for pre-commutative algebras, in our terminology. 

Namely, Definition \ref{defn:preP-algebra} applied to the operad $\Com $ governing the variety of associative 
and commutative algebras shows that $\pre\Com $ is defined by \eqref{eq:DendIdent} together with 
$a_1\succ a_2 = a_2\prec a_1$. Therefore, pre-commutative algebras may be considered as systems with 
one operation $a_1a_2 = a_1\succ a_2$ satisfying 
\begin{equation}\label{eq:Zinbiel}
 (a_1a_2 + a_2a_1)a_3 = a_1(a_2a_3).
\end{equation}

We will use the following notation: for $a_1,a_2,\dots , a_n\in A$, $A\in \pre\As$, 
let $[a_1\dots a_n]$ denote $(\dots ((a_1a_2)a_3)\dots a_n)$.

\begin{thm}\label{thm:ZinbielGSB}
The following polynomials in $\pre\As\langle X\rangle$ form a GSB of the free pre-commutative algebra $\pre\Com\langle X\rangle$:
\begin{gather}
 M_{u,v} =  u*v - uv - vu, 
  \label{eq:ZinbielGSB-1} \\
Z_{u,v,w} = u(vw) - (uv)w - (vu)w,
  \label{eq:ZinbielGSB-2} 
\end{gather}
where $u,v,w \in U$.   
\end{thm}

\begin{proof}
Denote the set of polynomials in the statement by~$S$.

Obviously, the compositions of $*$-inclusion (C1) are the compositions of
of 
$M_{u,v(wp)}$ and $Z_{v,w,p}$ (relative to $u*v(wp)$) or $M_{u(vw),p}$ and $Z_{u,v,w}$ (relative to $u(vw)*p$).
For example, the first of these compositions is equal to 
\[
h= u(v(wp))+(v(wp))u - u*((vw)p)-u*((wv)p) \to_S Z_{v,w,p}u + uZ_{v,w,p}\to_{\hat S} 0,
\]
so $h\equiv 0\pmod {\hat S, *, u*(v(wp))}.$
The second one is trivial by similar reasons. 

Let us compute a composition of $*$-intersection (C2): 
the composition of $M_{u,v}$ and $M_{v,w}$ (relative to $u*v*w$)
is equal to
\[
 h = (uv+vu)*w - u*(vw+wv).
\]
By definition, 
$h\to _{S,*} f = w(uv)+ (uv)w +w(vu)+(vu)w - u(vw) - u(wv) - (vw)u - (wv)u $.
The latter may be presented as 
\[
 f =  Z_{w,u,v} + Z_{w,v,u} - Z_{u,v,w} - Z_{u,w,v},
\]
so $h\equiv 0\pmod {\hat S ,*,u*v*w}$.

Consider compositions of type (C3).
Suppose $s_1 \in S$ and $s_2 \in S_0$ have a composition of $\succ $-inclusion. 
If $s_1\in S_0$ then one may apply Theorem~\ref{thm:CD_nonass} and the known result on the linear basis 
of the free Zinbiel algebra \cite{Loday01} to conclude that the set of $\{ Z_{u,v,w} \mid u,v,w\in X^{**}\}$
is a GSB in $M\langle X\rangle $, so it is closed with respect to all compositions of inclusion. 
However, it is not hard to compute all possible compositions straightforwardly and show they are all trivial. 

If $s_1 \in S_1$ then $s_1=M_{v,w}$ for $v,w\in U$. 
Assume $\hat s_2$ is a subword in $v$. Then the composition of $\succ$-inclusion is equal to 
\[
 h = vw+wv - (v_1+v_2)*w,
\]
where $v -v_1 - v_2 \in \hat S$, $v> v_1, v_2 $. 
Obviously, 
\[
 h \to_{S,*} vw+wv - v_1w-wv_1 - v_2w - wv_2 \equiv 0 \pmod {\hat S, *, v*w} 
\]
The second case, when $\bar s_2$ is a subword in $w$, is completely analogous.

Therefore, by Lemma~\ref{lem:C1--C3}  $\hat S$ 
is a GSB in $\As\langle U\rangle$. 
It remains to show that all compositions (C4) and (C5) of right and left multiplication belong to $I_*(\hat S)$.

Relation \eqref{eq:ZinbielGSB-1} shows that a composition of right multiplication of $M_{u,v}$
by $w$ is exactly $Z_{u,v,w}$.
For compositions of left multiplication, note that \eqref{eq:DendIdent-2} implies 
\begin{equation}\label{eq:LeftMul*-product}
  (w_1*\dots * w_k)(u*v) = [w_1,\dots ,w_k, u]*v -[w_1,\dots ,w_k, u]v + [w_1,\dots ,w_k, u,v], 
\end{equation}
in $\pre\As\langle X\rangle $, where $[x_1,\dots , x_m] = x_1(x_2 \dots (x_{m-1}x_m)\dots )$.
Therefore, the composition of left multiplication of $M_{u*v} $ by $w=w_1*\dots * w_k$ 
is equal to 
\[
h = [w_1,\dots ,w_k, u]*v -[w_1,\dots ,w_k, u]v - [w_1,\dots ,w_k, v,u].
\]
Obviously, $h\to_{S,*} f= [v ,w_1,\dots ,w_k, u] - [w_1,\dots ,w_k, v,u]$, 
where $f\in I_*(\hat S)$: for $k=1$
\[
f= [v,w_1,u]-[w_1,v,u] = Z_{v,w_1,u}- Z_{w_1,v,u}  ,
\]
and 
\[
f = Z_{v,w_1,[w_2,\dots , w_k, u]} - Z_{w_1,v,[w_2,\dots , w_k, u]} + w_1([v,w_2,\dots , w_k, u] - [w_2,\dots , w_k, u])\in I_*(\hat S)
\]
for $k>1$ by induction.

Hence, $S$ is a GSB in $\pre\As\langle X\rangle $, and monomials of the form 
$(\dots (x_1x_2)\dots x_n)$, $x_i\in X$, form a linear basis 
of $D=\pre\As\langle X\mid S\rangle $. Since all these monomials belong to $U$, 
the entire $D$ satisfies pre-commutativity (Zinbiel) identity, so 
$D\simeq \pre\Com\langle X\rangle $.
\end{proof}

\subsection{GSB of the universal enveloping pre-associative algebra of an associative algebra}

Every pre-associative algebra $P$ considered as a space with one product $x*y = x\prec y + x\succ y$ 
is an ordinary associative algebra denoted $P^{(+)}$.
For every associative algebra $A$ with a product denoted by $\cdot $ there exists its universal enveloping 
pre-associative algebra $U_*(A)$. 
Namely, $A$ maps into $U_*(A)^{(+)}$, $U_*(A)$ is generated (as a pre-algebra) by the image of $A$, and every 
homomorphism $A\to P^{(+)}$ of ordinary algebras may be extended to a homomorphism $U_*(A)\to P$ of pre-algebras.
It is clear that $U_*(A)$ may be presented by generators and relations as follows. Suppose $X$ is a linear basis of $A$, then 
\[
 U_*(A) \simeq  \pre\As \langle X\mid x*y - x\cdot y, x,y\in X\rangle.
\]

Introduce the following notation: 
$L_{w_1,\dots, w_k}(u) = [w_1,\dots , w_k, u]$ for $w_1,\dots, w_k, u\in U$, $k\ge 0$ 
(for $k=0$, denote the operator $L_{w_1,\dots, w_k}$ by $\mathbf 1$). 
It follows from 
\eqref{eq:LeftMul*-product} that for every $a,b\in X$ and for every $L=L_{w_1,\dots , w_k}$
\[
 M_{L,a,b} = L(a)*b -L(a)b + L(ab)-L(a\cdot b) 
\]
belongs to the ideal of $\pre\As\langle X\rangle $ generated by $x*y-x\cdot y$, $x,y\in X$. Note that an arbitrary word in $U$ may 
be uniquely presented as $L(a)$.
Moreover, 
\[
 W_{L,a,b,u} = M_{L,a,b}u = L(a)(bu)-(L(a)b)u + L(ab)u - L(a\cdot b)u
\]
also belongs to the same ideal. 

\begin{thm}\label{thm:GSB_of_preAsEnvelope}
The set of relations 
$M_{L,a,b}$, $W_{L,a,b,u}$, where $L=L_{w_1,\dots , w_k}$, $w_1,\dots , w_k\in U$ ($k\ge 0$),
$a,b\in X$, $u\in U$, 
is a GSB in $\pre\As\langle X\rangle $. 
\end{thm}

This is a GSB of the universal enveloping pre-associative algebra $U_*(A)$ of an associative algebra $A$.

\begin{proof}
Denote by $S$ the set of relations in the statement.
Let us compute all compositions (C1)--(C3) of $S$ to make sure $\hat S$ is a GSB 
in $\As\langle U\rangle $.
 
The only composition of $*$-intersection (C1) comes from $M_{\mathbf 1,a,b}$ and $M_{\mathbf 1, b,c}$, $a,b,c\in X$. 
It is trivial due to associativity of $A$.
 
Let us consider compositions (C2) and (C3) together since (C2) is a ``degenerate'' case of (C3).
Suppose $M_{L,a,b}$ and $W_{L',c,d,u}$ have a composition of $*$-inclusion or $\succ $-inclusion, 
$L=L_{w_1,\dots , w_k}$. There are two possible cases:
\begin{enumerate}
 \item $L'(c)(du)$ is a subword in a word $w_i$;
 \item $L'(c)(du)$ coincides 
with $L_{w_l,\dots, w_k}(a)$, $l\ge 1$ ($l=1$ corresponds to composition of $*$-inclusion).
\end{enumerate}

In the first case, there exist $w_{i1},w_{i2},w_{i3}\in U$ such that 
$w_i-w_{i1}-w_{i2}-w_{i3} \in \hat S_0$, $w_{ij}<w_i$. Then 
for $L_j = L_{w_1,\dots, w_{ij},\dots , w_k}$, $j=1,2,3$, we have 
$L(w)-L_1(w)-L_2(w)-L_3(w) \in \hat S_0$, $L_j(w)<L(w)$ for all $w\in U$.
The composition of $\succ $-inclusion has the form 
\[
 h = L(a)b - L(ab)+L(a\cdot b) -\sum\limits_{j=1}^3 L_j(a)*b,
\]
and, obviously, $h\equiv 0\pmod {\hat S , *, L(a)*b}$.

In the second case, $L(a)=[w_1,\dots, w_{l-1}, L'(c)(du)]$, so 
$u=L''(a)$, where $L''=L_{v_1,\dots , v_m}$ for some $v_1,\dots, v_m\in U$,
and $L(x) = [w_1,\dots, w_{l-1}, L'(c)(dL''(x))]$ for $x\in U$. 
Denote $L_{(1)}=L_{w_1,\dots, w_{l-1}}$.
Then the composition of $\succ$-inclusion (or $*$-inclusion for $l=1$) 
is equal to
 \begin{multline}\nonumber
  h = L(a)b-L(ab)+L(a\cdot b) - L_{(1)}\big( (L'(c)d)L''(a) - L'(cd)L''(a) + L'(c\cdot d)L''(a)  \big )* b \\
  =L(a)b-L(ab)+L(a\cdot b)-L_{(11)}(a)*b + L_{(12)}(a)*b -L_{(13)}(a)*b \\
  \equiv 
  L(a)b-L(ab)+L(a\cdot b)
  -L_{(11)}(a)b + L_{(11)}(ab)-L_{(11)}(a\cdot b ) \\
  + L_{(12)}(a)b - L_{(12)}(ab) + L_{(12)}(a\cdot b)   
  -L_{(13)}(a)b + L_{(13)}(ab)-L_{(13)}(a\cdot b )
  \pmod {S , * L(a)*b},
\end{multline}
where 
\[
\begin{gathered}
 L_{(11)} = L_{w_1,\dots, w_{l-1}, L'(c)d, v_1,\dots, v_m}, \\
 L_{(12)} = L_{w_1,\dots, w_{l-1}, L'(cd), v_1,\dots, v_m}, \\
 L_{(13)} = L_{w_1,\dots, w_{l-1}, L'(c\cdot d), v_1,\dots, v_m}.
\end{gathered}
\]
Note that
\[
 L(w) = L_{(1)}\big ( L'(c)(dL''(w))\big) 
 \equiv L_{(11)}(w) -L_{(12)}(w) + L_{(13)}(w) \pmod {\hat S, *, L(a)*b}
\]
for every $w\in U$. Hence, 
\[
 h\equiv 0 \pmod {\hat S, *, L(a)*b}.
\]

We have shown $\hat S$ to be a GSB in $\As\langle U\rangle$. To prove the theorem, it remains 
to check that all compositions of type (C4), (C5) of relations from $S$ are trivial modulo $\hat S$.

Obviously, the composition of right multiplication 
of $M_{L,a,b}$ by $u\in U$ is equal to $W_{L,a,b,u}$, i.e., it is trivial.
Moreover, 
$u M_{L,a,b} = M_{L',a,b}$, where
$L'(x) = uL(x)$, $u\in U$. Therefore, $S$ is a GSB in $\pre\As\langle X\rangle $.
\end{proof}

\begin{cor}
Linear basis of the universal enveloping pre-associative algebra $U(A)$ of an associative algebra $A$ 
consists of $u_1*u_2*\dots *u_n$, $n\ge 1$, $u_i\in U$, 
where $|u_2|,\dots , |u_n|>1$ and neither of $u_i$ contains a subword of the form $w(xv)$, $w,v\in U$, $x\in X$.
In particular, the pair of varieties $\As$, $\pre\As $ has PBW-property in the sense of \cite{MikhShest}.
\end{cor}

\subsection{Open problems}

One of the most intriguing problems related with application of GSB theory for pre-associative algebras
is related with the morphism of operads 
\[
({-}): \pre\Lie \to \pre\As
\]
defined by $xy\mapsto x\succ y - y\prec x$. Indeed, every pre-associative algebra $A$ with respect to new operation 
$a\cdot b = a\succ b -b\prec a$, $a,b\in A$, is a pre-Lie algebra $A^{(-)}$.

\medskip
\noindent
{\bf Problem 1.} 
Find the Gr\"obner---Shirshov basis of the universal pre-associative envelope 
\[
 U_{(-)}(L) = \pre\As \langle X \mid a*b - ab -ba + b\cdot a, a,b\in X \rangle, \quad \text{$X$ is a basis of $L$},
\]
of a pre-Lie algebra $L$.

\medskip
\noindent
{\bf Problem 2.}
Whether the triple $(\pre\As, \pre\Lie , ({-}))$ has the PBW-property?

\medskip
There is a morphism $\psi $ from the operad $\Brace$ governing the class of brace-algebras \cite{GersVoronov}
to the operad $\pre\As$ described in \cite{ChapBrace}. It was shown in \cite{ChapBrace} that 
the triple $(\pre\As, \Brace, \psi )$ has the PBW-property. 

\medskip
\noindent
{\bf Problem 3.}
Find the GSB of $U_\psi (B)$ for a brace-algebra $B$.

\end{document}